\documentclass{rspublic}

\input epsf.tex
\epsfverbosetrue

\begin{document}

\title[On Parrondo's Paradox]{
On Parrondo's paradox: how to construct unfair games by composing fair games}

\author[E.S.~Key, M.M.~K\l osek, D.~Abbott]
{Eric S.~Key$^1$, Ma\l gorzata M.~K\l osek$^1$, Derek Abbott$^2$}

\affiliation{$^1$Department of Mathematical Sciences, University of Wisconsin Milwaukee,
Milwaukee WI 53201 USA\\
$^2$Centre for Biomedical Engineering (CBME), Department of Electrical
\& Electronic Engineering, University of Adelaide, Adelaide SA 5005 Australia}

\label{firstpage}

\maketitle

\begin{abstract}{random walk in a periodic environment, random transport,
random games, Parrondo's paradox}
We construct games of chance from simpler games of chance. We show that it may
happen that the simpler games of chance are fair or unfavourable to a player and
yet the new combined game is favourable -- this is a counter-intuitive phenomenon
known as Parrondo's paradox.  We observe that  all of the games in question
are random walks in periodic environments (RWPE) when viewed on the proper time
scale. Consequently, we use RWPE techniques  to derive conditions under which
Parrondo's paradox occurs.
\end{abstract}

\section{Introduction}
Parrondian strategies are where losing games can cooperate to win
(Harmer \& Abbott 1999$a$). The original example of Parrondo's games consist
of two coin tossing games. Game A consists of Coin 1 biased to lose.
Game B consists of two coins -- Coin 2 with losing bias and Coin 3 with winning
bias -- but a state-dependent rule is chosen to favour the losing Coin 2. Hence
both games A and B are losing games. However when A and B are alternated in a
deterministic or even random manner, the player surprisingly has a winning
expectation. This effect has been interpreted in terms of a discrete-time
Brownian ratchet, at length, elsewhere (Harmer \& Abbott 1999$b$) -- where
conventional Brownian ratchets (Doering 1995) have been the inspiration. An
alternative view, we call the Boston interpretation (Stanley 1999, group
discussion), recognises that although game B favours Coin 2 with losing bias, if
the state-dependence is removed game B now favours the winning Coin 3 -- then
when games A and B are mixed, game A has the affect of randomisation or
`break-up' of game B's state-dependence, thus tilting favour towards Coin 3
with winning bias. This explanation was also independently deduced by
J.~Maynard Smith (1999, personal communication).
Further to the ratchet interpretation and Boston interpretation, this paper
will examine another viewpoint by considering the process as a random walk in a
periodic environment (RWPE).

We now briefly summarise the literature on Parrondo's games. In (Harmer
\textit{et al.} 2000$a$) the state-dependent rule, for game B, is to choose
Coin 2 if the player's capital is a multiple of some integer $M$ -- analysis
showed the paradox could hold for general values of  $M$. In (Pearce, 2000$a$)
it is shown that the paradox can hold when both games A and B are multiple
coin games. In (Pearce 2000$b$) we have the first analysis of Parrondo's games
in terms of Shannon entropy and in (Harmer \textit{et al.} 2000$b$) the entropy
parameter spaces are graphically displayed. The probability parameter space
is shown in (Harmer \textit{et al.} 2000$c$). In (Lee \textit{et al.}  2000) a
minimal three-state game with asymmetric transition probabilities is analysed
and in (Parrondo \textit{et al.} 2000) state-dependence on capital is
replaced by dependence on the past history of the game, leading to a larger
probability parameter space.

The surge of interest in analysing Parrondian games is motivated by a number
of areas. Information theorists have long studied the problem of producing
a fair game from biased coins (Gargamo and Vaccaro 1999) and the roots of
this  can be traced back to the work of von Neumann (1951) -- Parrondo's games
go a step further in producing a winning game from losing games. Seigman
(1999, personal communication) has reinterpreted `capital' of the games in terms
electron occupancies of energy levels -- the paradox can then be reproduced
using the rate equation approach typically used in laser analysis. In the
physical world there are many types of processes where losing helps to win, such
as a sacrifice in the game of chess or a valley in the fitness landscape of an
animal species. Many biological effects are linked to ratchet type phenomena
and Westerhoff {\it et al.} (1986) have analysed enzyme transport with a
four-state model. Applicability to population genetics, evolution and economics
has been suggested (McClintock 1999). In finance, Maslov and Zhang (1998) have
shown that under certain conditions capital can grow by investing in an asset
with {\it negative} typical growth rate. Quantum ratchets have now been
experimentally realised (Linke \textit{et al.} 1999) and recasting Parrondian
games, based on ratchet phenomena, as quantum games (Eisert \textit{et al.}
1999; Goldenberg \textit{et al.} 1999; Meyer 1999) is thus of interest.

In control theory, it can be shown that the combination of two unstable systems
can become stable (Allison \& Abbott 2000). Velocity of propagation through
an array of coupled oscillators, under certain conditions, can {\it increase}
even though the damping coefficient is {\it increased} (Sarmiento
\textit{et al.} 1999). In the area of granular flow, drift can occur in
a counter-intuitive direction such as exemplified in the famous Brazil nut
paradox (Rosato {\it et al.} 1987). Also declining branching processes can be combined to
increase (Key 1987). Plaskota (1996) shows that noisy information can sometimes
be better than clean information. In (Pinsky and Scheutzow 1992) it is shown
that with switched diffusion processes in random media it is possible to get a
positive-recurrent processes (\textit{i.e.}~with no drift) from mixed transient processes
(\textit{i.e.}~with drifts all in the same direction) -- this is almost certainly a
continuous time analogue to the Parrondian discrete-time process. Assuming we
construct Parrondo's games to only deal in transactions of one unit of capital
per event, then we have a skip-free process, and a statistical interpretation of
the central result is that declining birth-death processes can be combined to
form an increase.

In this paper we further investigate Parrondo's paradox. We construct a
class of composite games and investigate their fairness by formulating the
problem in the language of random walks in periodic environments (RWPE). We
find many new, interesting, and counter-intuitive results.

\section{Mathematical tools}

We construct a composite game from two simple games $A$ and $B$.
These two games can be combined in two ways: deterministically or stochastically.

If we have played $n$ times,
and $n$ is divisible by $k$ (for an integer $k$), we then play game $A$.
If we have played $n$
times and $n$ is not divisible by $k$, we play game $B$.  Thus games $A$ and
$B$ are alternated in a deterministic pattern.  We denote by $Y_n$ our capital
after $n$ plays of this game.

To alternate games $A$ and $B$ randomly we toss a coin with probability $p$ of
heads.  If the coin comes up heads, we play game $A$, and if tails, we play
game $B$.
We denote by $Z_n$ be our capital after $n$ repetitions.

The sequence of values of our capital in either of these games
is a Markov random walk which changes by $\pm 1$ in each epoch.
The two random walks $Z_n$ and $Y_n$ differ
in that $Z_n$ is time homogeneous, while $Y_n$ is not time homogeneous.
Moreover, $Z_n$ is a random walk in a periodic environment. The process
 $Y_n$ is not a RWPE, but the process $Y'_n \equiv Y_{kn}$ is a RWPE.

As shown in \S4 we extend this construction to define composite games from more
than two simple games. In each case we identify the capital of the player as a
RWPE.

We say that a game is fair, winning or losing if the random walk for the capital of a
player, $X_n$, is recurrent or transient to $\infty$, or to $-\infty$,
respectively. That is,
a Markov chain $X_n$ is

\begin{center}
\begin{tabular}{l}
recurrent (fair) if $\displaystyle P\{-\infty = \lim\inf_{n \rightarrow
\infty}X_n < \lim\sup_{n \rightarrow \infty}X_n = \infty$\} = 1; \\
transient to $\infty$ (winning) if $\displaystyle P\{\lim_{n \rightarrow\infty}X_n = +\infty\} =
1$; \\
transient to -$\infty$ (losing) if $\displaystyle P\{\lim_{n \rightarrow\infty}X_n =
-\infty\}$=1.
\end{tabular}
\end{center}

We note that the characterisation of a game as fair, winning, or
losing  by the traditional
comparisons $E[X_{n+1} | X_n] = X_n$, $E[X_{n+1} | X_n] > X_n$, and $E[X_{n+1}
| X_n] < X_n$, respectively, does not cover the behaviour
all random walks, in particular RWPE's.

\subsection{Key's criterion}

We consider a time homogeneous random walk, $X_n$, in an $N$--periodic
environment, or equivalently we have a state dependent random game ${\bf W}$.
We assume that the maximal step size in the positive direction is $R$,
while in the negative direction it is $L$, that is $P\{X_{n+1} \in \{-L + k, \ldots, x + R\} | X_n =
k\} = 1$. Moreover, for each $k$ the maximum right and left step sizes are always possible,
that is, $P\{X_{n+1} =
-L + k | X_n = k\}P\{X_{n+1} = R + k | X_n = k\} > 0$. Given the environment, the walk
$X_n$ obeys the backward master equation

\begin{eqnarray}
P\{\cdot\, | X_n = k\} = \sum_{j =  -L}^{R}e(k, j)P\{\cdot\, | X_{n + 1} = j + k\},
\label{BME}
\end{eqnarray}

where $e(k, j) \equiv P\{X_{n+1} = j + k | X_n = k)$ denotes the transition probability from the
state $k$ to $j + k$ in one time epoch.
We denote by $f_{k - i} \equiv P\{\cdot\, | X_{n + 1} = k - i\}$ (since $X_n$ is time
homogeneous $f$ does not depend on $n$) and rewrite equation (\ref{BME}) as a system
for the vector $[f_{-L + k}, f_{-L + 1 + k}, \ldots, f_{R - 1 + k}]^T$ with the
matrix ${\bf A}_k$ whose entries given by

\begin{eqnarray}
{\bf A}_k[i, j] = \left\{
\begin{array}{ll}
-e(k, -L+j)/e(k, -L)&\mbox{if\ }i=1, j\neq L\\
(1-e(k, 0))/e(k, -L)&\mbox{if\ }i=1, j=L\\
1& \mbox{if\ } i \geq 2, j = i-1\\
0&\mbox{otherwise}
\end{array}\right.. \nonumber
\end{eqnarray}

We also define the matrix ${\bf M} = {\bf A}_1{\bf
A}_{2}\ldots{\bf A}_{n}$.

According to (Key 1984) we define constants $d_i, i= 1,2,
\ldots R+L$ as follows. For each eigenvalue $\lambda_i$ of ${\bf M}$ (including
multiplicities), we put $d_i = \log(|\lambda_i|)$, and we list the $d_i$'s in
increasing order, so that $d_1 \leq d_2 \leq \ldots \leq d_{R + L}$. Then
\begin{eqnarray}
\hspace*{-14ex}&&\hspace*{-6ex}\mbox{if\ } \ln(c({\bf W})) \equiv d_R + d_{R + 1} > 0
\mbox{\ then\ the\ RWPE\ } X_n
\mbox{\ is transient to\ } \infty ; \nonumber\\
&&\hspace*{-6ex}\mbox{if\ } \ln(c({\bf W})) \equiv d_R + d_{R + 1} = 0
\mbox{\ then\ the\ RWPE\ } X_n
\mbox{\ is recurrent\ } ;\nonumber\\
&&\hspace*{-6ex}\mbox{if\ } \ln(c({\bf W})) \equiv d_R + d_{R + 1} < 0
\mbox{\ then\ the\ RWPE\ } X_n
\mbox{\ is transient to\ } -\infty. \label{CRIT}
\end{eqnarray}

It is shown in (Key\&K\l osek 2000) that $X_n$ is recurrent if the
characteristic polynomial of ${\bf M}$ has a double root at $1$.

\section{ Games composed of 2-periodic games}

First we consider an example of two 2-periodic games ${\bf P} = (P_0, P_1)$ and
${\bf Q} = (Q_0, Q_1)$; that is, we have two RWPE's, each of them
2-periodic (in space). If $X_n$ denotes the capital of a player at time $n$ which plays the game
according to the rule ${\bf P}$, then the transition probabilities are given by

\begin{eqnarray}
&&P\{X_{n+1} = x + 1 | X_n = x\} = \left\{
\begin{array}{ll}
p_0& \mbox{if } x = 0 \mbox{\ mod\ }2\\
p_1&\mbox{if } x = 1\mbox{\ mod\ } 2
\end{array}\right.\nonumber\\
&&P\{X_{n+1} = x - 1 | X_n = x\} = \left\{
\begin{array}{ll}
1 - p_0& \mbox{if } x = 0 \mbox{\ mod\ } 2\\
1 - p_1&\mbox{if } x = 1\mbox{\ mod\ } 2
\end{array}\right. ; \label{TRP2P}
\end{eqnarray}

that is, if the capital is even it changes according to $P_0$, and it is odd it
changes according to $P_1$ -- cf.~figure~1.  The transition probabilities for the RWPE governed
by ${\bf Q}$ are given by formulas analogous to (\ref{TRP2P}), with $p's$
replaced by $q's$, with the rule $Q_0$ and $Q_1$ at the even and odd
positions, respectively.

\begin{figure}
\epsfxsize=4in\epsfysize=1.5in\epsffile{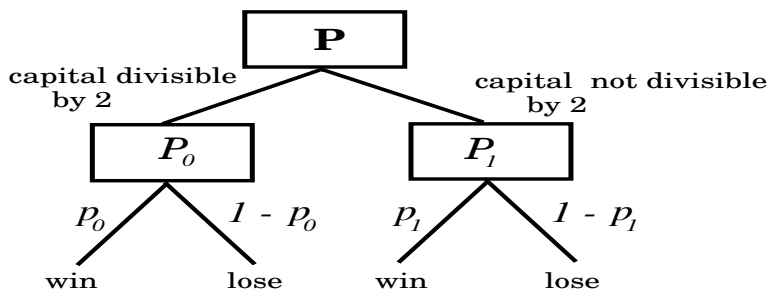}
\caption{Schematic representation of a 2--(state) periodic game ${\bf P} = (P_0, P_1)$.
Games composed of such games are discussed in \S 3.}
\end{figure}

To employ the criterion of (Key 1984) to the RWPE governed by ${\bf P}$ we
construct the matrix ${\bf M}$

\begin{eqnarray}
{\bf M} = {\bf A}_1{\bf A}_0 \equiv \left[
\begin{array}{cc}
\displaystyle 1/(1 - p_1)& -p_1/(1 - p_1)\\
1&0
\end{array}\right]\left[
\begin{array}{cc}
\displaystyle 1/(1 - p_0)& -p_0/(1 - p_0)\\
1&0
\end{array}\right]. \nonumber
\end{eqnarray}

Since here $R = L = 1$, $d_R + d_{R+1} = d_1 + d_2 = \log(|\det({\bf M})|)$,
so we find that the RWPE
governed by ${\bf P}$ is transient to $-\infty$, recurrent, or transient to
$+\infty$ if
$\displaystyle c({\bf P}) \equiv c(P_0, P_1) \equiv
\frac{p_0p_1}{(1 - p_0)(1 - p_1)} < (=) (>) 1. $ We note that this analysis easily extends
to any period $N$ of the environment -- cf. \S4.

\subsection{An example of a pseudo--paradox}
We suppose that the games ${\bf P}$ and ${\bf Q}$ are fair, that is
$\displaystyle c(P_0, P_1) = p_0(1 - p_0)/(p_1(1 - p_1)) = 1$ and
$c(Q_0, Q_1) = q_0(1 - q_0)/(q_1(1 - q_1)) = 1$. If we
alternate (deterministically) games ${\bf P}$ and ${\bf Q}$ in that very order,
then, starting at the origin, the composite game may be fair (losing) (winning)
if $c({\bf P}, {\bf Q}) = p_0q_1/((1 - p_0)(1 - q_1))
= ( < ) ( > ) 1$. Under this strategy we play the game $P_0$ when the winnings
are even and the game $Q_1$ when the winnings are odd.  The games $P_1$ and
$Q_0$ are never played. Hence $c({\bf{P}}, {\bf{Q}}) = c(P_0, Q_1)$, and the
fact that the composite game is losing when it is constructed from two winning
games is just an apparent paradox. If we play the same two games when starting
at an odd position then $c({\bf P}, {\bf Q}) = c(P_1, Q_0) = 1/c(P_0, Q_1)$,
and we have a winning game constructed out of two winning games.

\subsection{Effects of randomisation}

\begin{figure}
\epsfxsize=4in\epsfysize=3.0in\epsffile{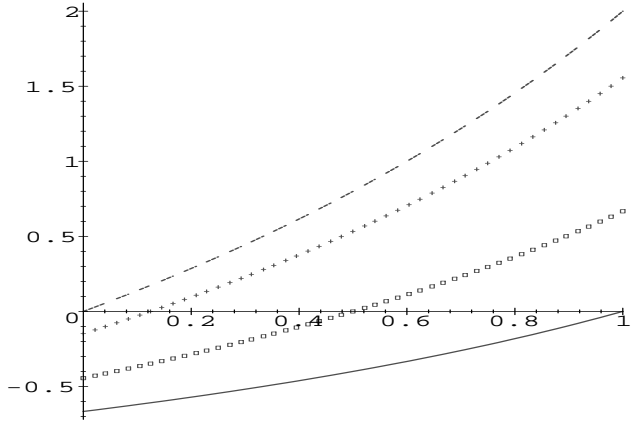}
\caption{We take two 2--periodic fair games with $p_1 = 1/2$ and $q_1 = 1/4$, and plot $
c({\bf P}, {\bf Q}) - 1$ in Eq.(\ref{C2G2PR}) as a function of $g_1$ for various values of
$g_0$; that is, the composite game is winning for positive values on the graph, and losing
for negative values. We have: dashes: $g_0 = 0$; crosses: $g_0 = 1/8$; boxes: $g_0 = 1/2$;
solid: $g_0 = 1$.}
\end{figure}

Now we construct a RWPE, $X_n$, by choosing at random between ${\bf P}$
and ${\bf Q}$.
If our capital is even then a coin with probability of heads
$g_0$ is used to choose whether the game $P_0$ or $Q_0$ will be played; if our
capital is odd then a coin with probability of heads $g_1$ is
used to choose between the games $P_1$ and $Q_1$. In this scheme all four games
are played. The transition probabilities are given by

\begin{eqnarray}
P\{X_{n + 1} = x + 1 | X_n = x\} = \left\{
\begin{array}{ll}
g_0p_0 + ( 1- g_0)q_0& \mbox{if } x = 0 \mbox{\ mod\ } 2\\
g_1p_1 + (1 - g_1)q_1& \mbox{if } x = 1 \mbox{\ mod\ } 2
\end{array}\right. .\label{R2P2G}
\end{eqnarray}

The fairness of the composite game is determined by the factor
$c({\bf P}, {\bf Q})$ given by

\begin{eqnarray}
c({\bf P}, {\bf Q}) = \frac{g_0p_0 + ( 1- g_0)q_0}{(1 - g_0p_0 - ( 1-
g_0)q_0)}\frac{g_1p_1 + (1 - g_1)q_1}{(1 - g_1p_1 - (1 - g_1)q_1)}.
\label{C2G2PR}
\end{eqnarray}

By direct simplification of equation(\ref{C2G2PR}) we observe that if $g_0 = g_1$
(and the games ${\bf P}$ and ${\bf Q}$ are fair)
then the composite game is fair. Also, if $p_1 = q_1$ then the composite game
is fair. However, for other values of parameters two fair games can be used to compose a
game which is winning, losing or fair depending on how the two simple games are
randomised. Figure 2 illustrates this statement.
That is, randomisation can produce Parrondo's
paradox.

We also consider a game composed from two unfair games ${\bf P}$ and ${\bf Q}$. Games
${\bf P}$ and ${\bf Q}$ are taken to be unfair the same way; that is both are losing or
both are winning.
We note that when
one coin is used to choose whether ${\bf P}$ or ${\bf Q}$ is played, 
(\textit{i.e.},~when $g_0 = g_1$ in
equation (\ref{R2P2G})) then the composed game is always unfair. However, if two coins
are used then the randomised game may be winning when both simple games ${\bf P}$ and
${\bf Q}$ are losing. This
example of Parrondo's paradox is illustrated in figure 3.

\begin{figure}
\epsfxsize=4in\epsfysize=3.3in\epsffile{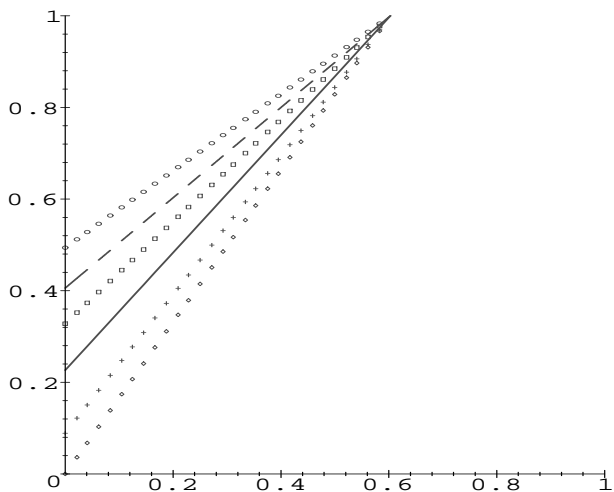}
\caption{We take two 2--periodic losing games with $p_0 = 0.675$, $p_1 = 0.1$, $q_1 = 0.75$ and
plot the fairness curve $c({\bf P}, {\bf Q}) = 1$ in Eq.(\ref{C2G2PR}) in the $(g_1,
g_0)$ plane for various values of $q_0$: circles: $q_0 = 0$; dashes: $q_0 = 0.075$;
 boxes: $q_0 = 0.125$; solid: $q_0 = 0.175$; crosses: $q_0 = 0.225$;
diamonds: $q_0 = 0.25$. The composite game is winning  when $g_0$ and $g_1$ are selected
above the fairness curve.}
\end{figure}

\subsection{Three 2-periodic games}

We consider three 2--periodic games ${\bf P}, {\bf Q}$ and ${\bf R}$, each
defined analogously to (\ref{TRP2P}), with transition probabilities given in
terms of $p_0, p_1$, $q_0, q_1$, and $r_0, r_1$, respectively. We construct the
composite game ${\bf PQR}$. We investigate whether the composite game may be
winning (losing) if the individual games ${\bf P}, {\bf Q}$ and ${\bf R}$ are
fair.  If the capital of the player at time $n$ is $X_n$ then the capital at
times $3n$, $Y_{n} \equiv X_{3n}$, is a random walk in a 2--periodic
environment, taking steps of  $\pm 3, \pm 1$. Specifically, we have

\begin{eqnarray}
&&\hspace*{-3ex}P\{Y_{n + 1} = x + 3 | Y_n = x\} = \left\{
\begin{array}{cc}
a_0 \equiv p_0q_1r_0  & \mbox{if\ }x = 0 \mbox{\ mod\ } 2\\
a_1 \equiv  p_1q_0r_1 &\mbox{if\ }x = 1 \mbox{\ mod\ } 2
\end{array}\right.\nonumber\\
&&\hspace*{-3ex}P\{Y_{n + 1} = x - 3 | Y_n = x\} = \left\{
\begin{array}{cc}
b_0  \equiv (1 - p_0)(1 - q_1)(1 - r_0) & \mbox{if\ }x = 0 \mbox{\ mod\ } 2\\
b_1 \equiv (1 - p_1)(1 - q_0)(1 - r_1) &\mbox{if\ }x = 1\mbox{\ mod\ } 2
\end{array}\right.\nonumber\\
&&\hspace*{-3ex}P\{Y_{n + 1} = x + 1 |  Y_n = x\} = \nonumber\\
&&\qquad\left\{
\begin{array}{cc}
c_0 \equiv p_0q_1(1 - r_0) + p_0(1 - q_1)r_0 + (1 - p_0)q_1r_0 & \mbox{if\ }x = 0 \mbox{\ mod\ } 2\\
c_1 \equiv p_1q_0(1 - r_1) + p_1(1 - q_0)r_1 + (1 - p_1)q_0r_1  &\mbox{if\ }x = 1\mbox{\ mod\ } 2
\end{array}\right.\nonumber\\
&&\hspace*{-3ex}P\{Y_{n + 1} = x - 1 |  Y_n = x\} =
\nonumber\\
&&\hspace*{-4ex}\left\{
\begin{array}{ll}
\hspace*{-1ex}d_0 \equiv (1 - p_0)(1 - q_1)r_0 + (1 - p_0)q_1(1 - r_0) + p_0(1 - q_1)(1 - r_0)
&\hspace*{-1ex} x = 0
\mbox{\ mod\ } 2\\
\hspace*{-1ex}d_1 \equiv (1 - p_1)(1 - q_0)r_1 + (1 - p_1)q_0(1 - r_1) + p_1(1 - q_0)(1 - r_1)
&\hspace*{-1ex}x = 1\mbox{\ mod\ }
2.\end{array}\right. \nonumber
\end{eqnarray}

Moreover, the random walk $Z_n \equiv Y_{2n}$ is an ordinary random walk, that
is a sum of iid random variables taking
values $\pm 6, \pm 4, \pm 2$, and 0, and the transition probabilities are

\begin{eqnarray}
&&P\{Z_{n + 1} = z + 6 | Z_n = z\} = a_0a_1\nonumber\\
&&P\{Z_{n + 1} = z - 6 | Z_n = z\} = b_0b_1\nonumber\\
&&P\{Z_{n + 1} = z + 4 | Z_n = z\} = a_0c_1 + c_0a_1\nonumber\\
&&P\{Z_{n + 1} = z - 4 | Z_n = z\} = b_0d_1 + d_0b_1\nonumber\\
&&P\{Z_{n + 1} = z + 2 | Z_n = z\} = c_0c_1 + a_0d_1 + a_1d_0\nonumber\\
&&P\{Z_{n + 1} = z - 2 | Z_n = z\} = d_0d_1 + b_0c_1 + b_1c_0\nonumber\\
&&P\{Z_{n + 1} = z |  Z_n = z\} = d_0c_1 + d_1c_0 + a_0b_1 + a_1b_0 . \nonumber
\end{eqnarray}

If the games ${\bf P}, {\bf Q}$ and ${\bf R}$ are fair, then by direct
calculations, $EZ_n = 0$, so the composite game is fair, and there is no
paradox in this case.

\section{Games composed of 3-periodic games}

We consider games composed of 3--periodic games ${\bf P}, {\bf Q}$, and ${\bf
R}$.  We define a 3--periodic game ${\bf P} = (P_0, P_1, P_2)$ by its
transition probabilities as

\begin{eqnarray}
P\{X_{n + 1} = x + 1 | X_n = x \} = \left\{
\begin{array}{ll}
p_0 & \mbox{\ if\ } x = 0 \mbox{\ mod \ } 3\\
p_1 & \mbox{\ if\ } x = 1 \mbox{\ mod \ } 3\\
p_2 & \mbox{\ if\ } x = 2 \mbox{\ mod \ } 3
\end{array}\right. ,\label{TRP3PW}
\end{eqnarray}
compare figure~4. Transition probabilities of the games ${\bf Q}$ and ${\bf R}$ are defined by
formulas analogous to (\ref{TRP3PW}) with $p'$s replaced by $q$'s and $r$'s,
respectively.  According to (Key 1984), to determine the fairness of the game
${\bf P}$ we consider the product of the eigenvalues of the matrix ${\bf M} \equiv {\bf A}_1{\bf
A}_2{\bf A}_0$ where

\begin{eqnarray}
{\bf A}_i =
\left[
\begin{array}{cc}
\displaystyle 1/(1 - p_i)& -p_i/(1 - p_i)\\
1&0
\end{array}\right] \quad i = 0, 1, 2.\nonumber
\end{eqnarray}

\begin{figure}
\epsfxsize=4in\epsfysize=1.5in\epsffile{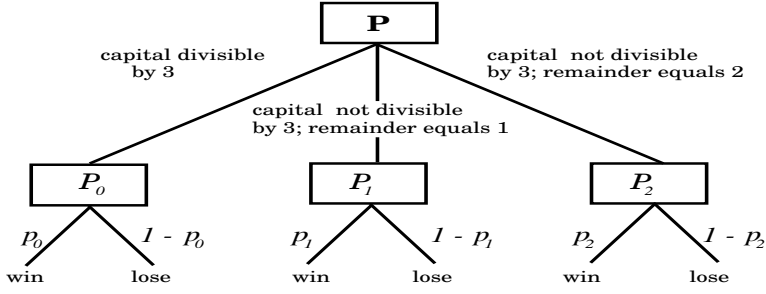}
\caption{Schematic representation of a 3--(state) periodic game ${\bf P} = (P_0,
P_1, P_2)$.
Games composed of such games are discussed in \S 4.}
\end{figure}

If
$\displaystyle{\frac{p_0p_1p_2}{(1 - p_0)(1 - p_1)(1 - p_2)} = ( < ) ( > ) 1}$
then the game ${\bf P}$ is fair (losing) (winning).

\subsection{Two 3--periodic games}

To analyse a game $({\bf P}, {\bf Q})$, where ${\bf P}$ and ${\bf Q}$ are
3--periodic games, we construct a random walk $Y_n \equiv X_{2n}$ where $X_n$
denotes the position of the walker at time $n$. The transition probabilities
for $Y_n$ are given by

\begin{eqnarray}
&&\hspace*{-3ex}P\{Y_{n + 1} = x + 2 | Y_n = x\} = \left\{
\begin{array}{ll}
 a_0 \equiv p_0q_1 &\mbox{\ if\ } x = 0 \mbox{\ mod 3}\\
a_1 \equiv p_1q_2  & \mbox{\ if\ } x = 1\mbox{\ mod 3}\\
a_2 \equiv p_2q_0  &\mbox{\ if\ } x = 2\mbox{\ mod 3}
\end{array}
\right.\nonumber\\*[2ex]
&&\hspace*{-3ex}P\{Y_{n + 1} = x | Y_n = x\} = \left\{
\begin{array}{ll}
b_0 \equiv p_0(1 - q_1) + (1 - p_0)q_2  &\mbox{\ if\ } x = 0 \mbox{\ mod 3}\\
b_1 \equiv p_1(1 - q_2) + (1 - p_1)q_0  &\mbox{\ if\ } x = 1\mbox{\ mod 3}\\
b_2 \equiv p_2(1 - q_0) + (1 - p_2)q_1  &\mbox{\ if\ } x = 2\mbox{\ mod 3}
 \end{array}
\right.\nonumber\\*[2ex]
&&\hspace*{-3ex}P\{Y_{n + 1} = x - 2 | Y_n = x\} = \left\{
\begin{array}{ll}
c_0 \equiv  (1 - p_0)(1 - q_2) &\mbox{\ if\ } x = 0 \mbox{\ mod 3}\\
c_1 \equiv (1 - p_1(1 - q_0) &\mbox{\ if\ } x = 1\mbox{\ mod 3}\\
c_2 \equiv (1 - p_2)(1 - q_1) &\mbox{\ if\ } x = 2\mbox{\ mod 3}.
\end{array}
 \right.\nonumber
\end{eqnarray}

Hence, $Y_n$ is a 3-periodic random walk taking steps $\pm 2$ and zero. To
employ Key's criterion most efficiently we observe that it is sufficient to
analyse $Y_n$ as a RWPE on the even integers which visits only nearest neighbours. To this end
we construct the $2\times 2$ matrix ${\bf M} = {\bf
A}_1{\bf A}_2{\bf A}_0$ where
\begin{eqnarray}
{\bf A}_i = \left[\begin{array}{cc}
(1 - b_i)/c_i&-a_i/c_i\\
1&0
\end{array}\right] \quad i = 0, 1, 2.
\nonumber
\end{eqnarray}

The determinant of ${\bf M}$ is given by
 $\det({\bf M})  = p_0p_1p_2q_0q_1q_2/[(1 - p_0)(1 - p_1)(1 - p_2)(1 -
q_0)(1 - q_1)(1 - q_2)]$. Hence, if the games ${\bf P}$ and ${\bf Q}$ are fair,
the composite game is fair and no paradox is observed.

We note that this reduction in the dimension of ${\bf M}$ occurs whenever the temporal
period is even, since the step sizes of the derived process $Y_n$ are even.
Specifically, if the temporal period is $T$ is even then ${\bf M}$ can be taken
to be the product of $T\times T$ matrices, and when the temporal period $T$ is
odd, then $\bf M$ is the product of $2T\times 2T$ matrices.

\subsection{Two 3--periodic games randomised}

\begin{figure}
\epsfxsize=4in\epsfysize=3.0in\epsffile{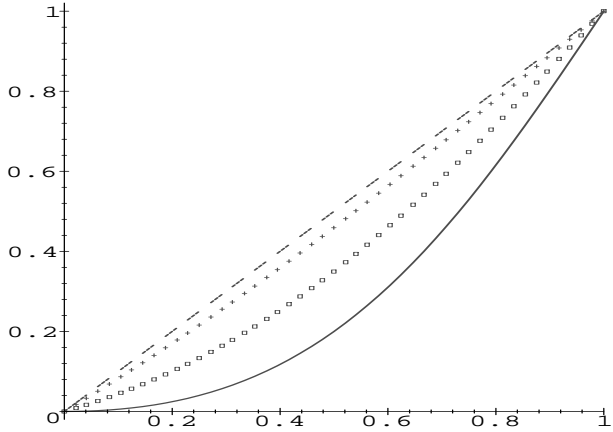}
\caption{We take two fair 3--periodic games with $g_1 =
g_2$, $p_1 = p_2 = 1/2$, and $q_1 = q_2$. For a fixed value of $q_1$ we plot
curves along which the composite game is fair, in the space $(g_1, g_0)$. As
$q_1$ varies from 0 to 1/2 the composite game is losing above and winning below
its corresponding ``fairness'' curve; solid: $q_1$ = 0; boxes: $q_1 = 1/6$;
crosses: $q_1 = 1/3$; dashes: the line $g_0 = g_1$.  As $q_1$ varies from 1/2
to 1 the composite game is winning below and losing above its corresponding
``fairness'' curve. In this case the problem is symmetric, so that we have:
solid: $q_1 = 1$; boxes $q_1 = 5/6$; crosses: $q_1 = 2/3$.}
\end{figure}

\begin{figure}
\epsfxsize=4in\epsfysize=3in\epsffile{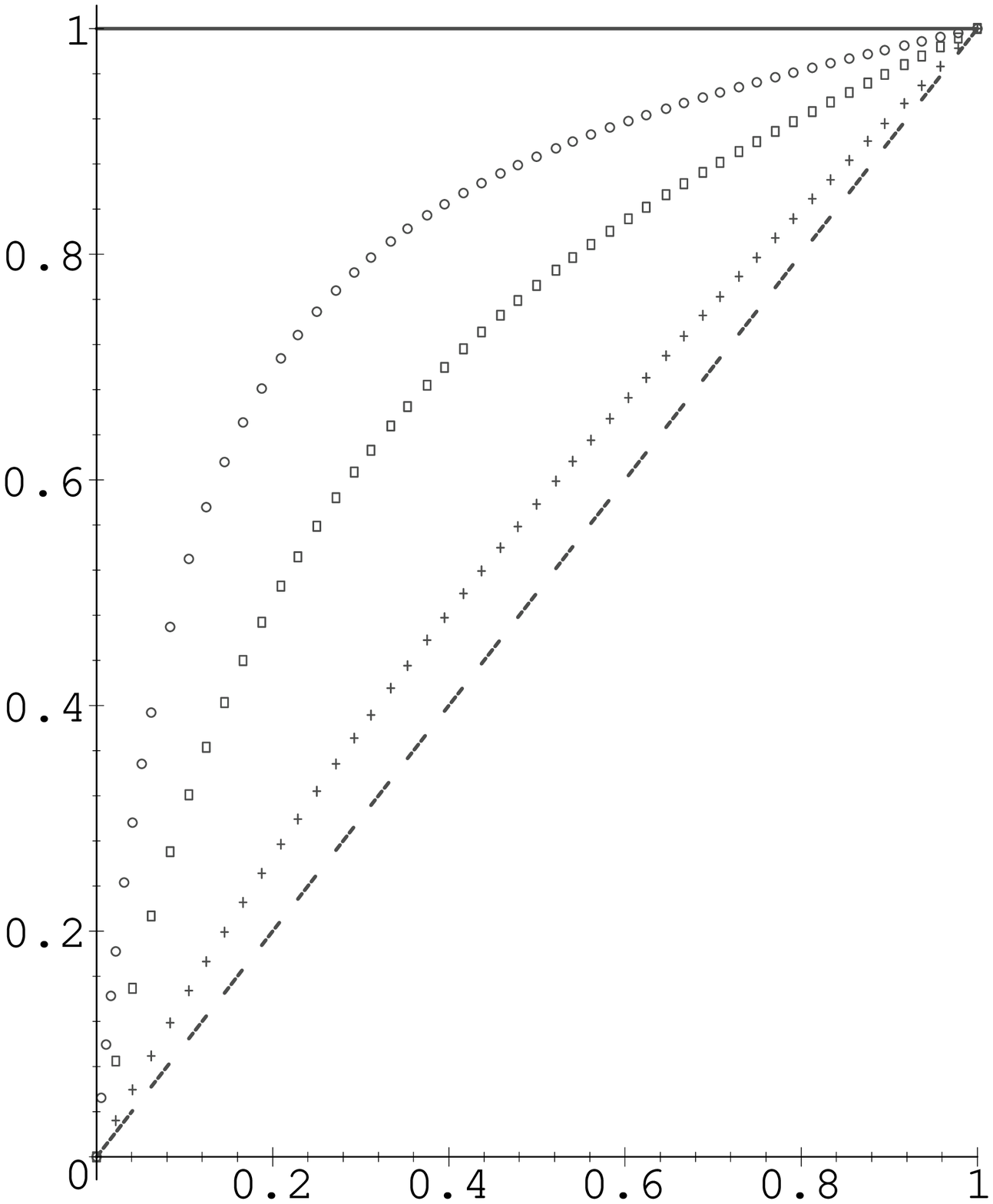}
\caption{We take two fair 3--periodic games with $g_1 =
g_2$, $p_1 = p_2 = 1/2$, and $q_0 = q_2$. We plot the curves along which the
composite game is fair in the space $(g_1, g_0)$, for a fixed value of
$q_0$. As $q_0$ varies from 0 to 1/2, the game is winning above and losing
below each of its ``fairness'' curves; solid: $q_0 = 0$; circles: $q_0 = 1/12$;
boxes: $q_0 = 1/6$; crosses: $q_0 = 1/3$; dashes: the line $g_0 = g_1$.  As
$q_0$ increases from 1/2 to 1 the game is losing above and winning below its
fairness curve. The case of this example is symmetric so that we have: crosses:
$q_0 = 2/3$; boxes: $q_0 = 5/6$; circles: $q_0 = 11/12$.}
\end{figure}

We construct a composite game from two 3--periodic random games ${\bf P}$ =
($P_0$, $P_1$, $P_2)$ and ${\bf Q} = (Q_0, Q_1, Q_2)$ by selecting at random the
game to be played at each step, with probability of playing $P_i$ equal to
$g_i, i = 0, 1, 2$. That is, we define a 3--periodic random walk $Y_n$ taking
values $\pm 1$ with transition probabilities

\begin{eqnarray}
P\{Y_{n + 1} = x + 1 | Y_n = x\} = \left\{\begin{array}{ll}
\rho_0 \equiv g_0p_0 + (1 - g_0)q_0  &\mbox{\ if\ } x = 0 \mbox{\ mod 3}\\
\rho_1 \equiv g_1p_1 + (1 - g_1)q_1  &\mbox{\ if\ } x = 1\mbox{\ mod 3}\\
\rho_2 \equiv g_2p_2 + (1 - g_2)q_2  &\mbox{\ if\ } x = 2\mbox{\ mod 3}
\end{array}\right. .\label{TP3PR}
\end{eqnarray}

The walk $Y_n$ is recurrent (transient to $\infty$) (transient to $-\infty$) according to
$\displaystyle c({\bf P}, {\bf Q})  \equiv \frac{\rho_0\rho_1\rho_2}
{(1 - \rho_0)(1 - \rho_1)(1 - \rho_2)}
 = ( > ) ( < ) 1$.
We investigate the value of $c({\bf P}, {\bf Q})$ for various values of the
parameters. We note that if the games ${\bf P}$ and ${\bf Q}$ are fair with
$p_1 = p_2 = q_1 = q_2$ then the randomised composite game is fair for all
possible choices of $g_i$'s. Next we take two fair games with fixed values $p_1
= p_2$ and $q_1 = q_2$,
that is ${\bf P} = (P_0, P_1,
P_1)$ and ${\bf Q} = (Q_0, Q_1, Q_1)$.

We investigate how randomisation affects the
fairness of the composite game. As shown in figure 5 there is a wide range of parameters
for which the randomised game is unfair, so that randomisation induces the paradox. We
note that the range of the randomisation parameters $g's$ is affected by the way games
${\bf P}$ and ${\bf Q}$ are structured. In the examples illustrated in figure 6 we take two
fair games ${\bf P} = (P_0, P_1, P_1)$ and ${\bf Q} = (Q_0, Q_1, Q_0)$ and we note that
the range of parameters $g_0, g_1 = g_2$ for which the randomised game composed from
these ``shifted'' games is larger compared to the previous ``not--shifted'' case --
cf. figures 5 and 6.

\subsection{Three 3--periodic games}

\begin{figure}
\epsfxsize=4in\epsfysize=3in\epsffile{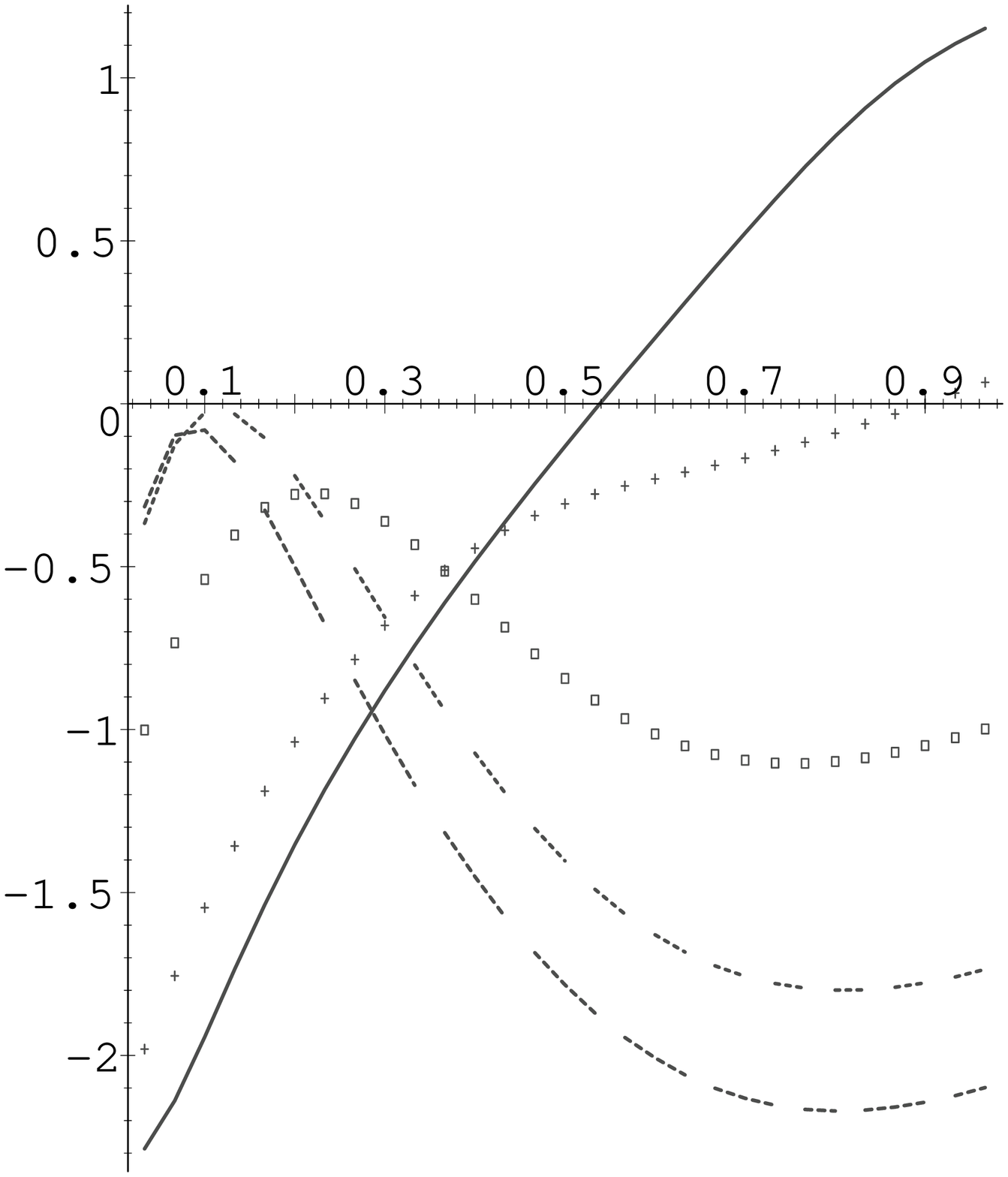}
\caption{We take 3 fair games with $p_0 = p_1 \equiv
p$, $q_0 = q_1 \equiv q$, and $r_0 = r_1 \equiv r$, and plot $\ln(c({\bf P}, {\bf
Q}, {\bf R})) = d_3 + d_4$ -- compare Eq.(\ref{CRIT}), as a function of $p$,
for $q = 0.1$ and various values of $r$, solid: $r = 3/4$; crosses: $r = 1/2$,
boxes: $r = 1/4$; short dashes: $r = 1/8$; long dashes: $ r = 1/16$. We note
that composite game is always losing, no matter what $p$ is, for a wide range
of values of the parameter $r$.  If $\ln(c({\bf P}, {\bf Q}, {\bf R}))$ is
positive (zero) (negative) then the game $({\bf P},{\bf Q}, {\bf R}))$ is
winning (fair) (losing).}
\end{figure}

\begin{figure}
\epsfxsize=4in\epsfysize=3in\epsffile{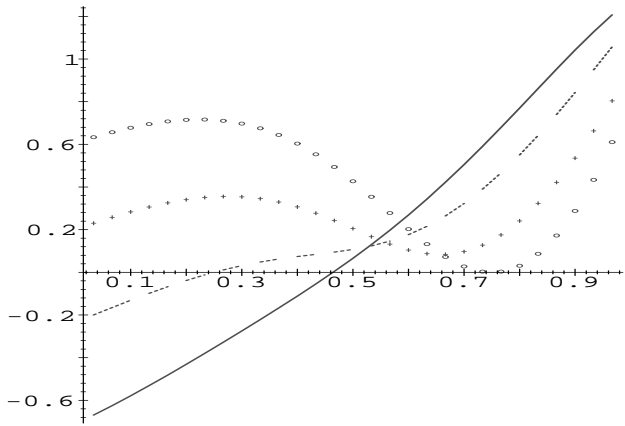}
\caption{We take 3 fair games with $p_0 = p_1 \equiv
p$, $q_0 = q_1 \equiv q$, and $r_0 = r_1 \equiv r$, and plot $\ln(c({\bf P}, {\bf
Q}, {\bf R})) = d_3 + d_4$ -- compare (\ref{CRIT}), as a function of $p$, for
$q = 3/4$ and various values of $r$: circles: $r = 3/4$, crosses: $r = 5/8$; dashes: $r
= 1/2$; solid: $r = 3/8$.
If $\ln(c({\bf P},
{\bf Q}, {\bf R}))$ is positive (zero) (negative) then the game $({\bf P},{\bf
Q}, {\bf R}))$ is winning (fair) (losing).}
\end{figure}

\begin{figure}
\epsfxsize=4in\epsfysize=3in\epsffile{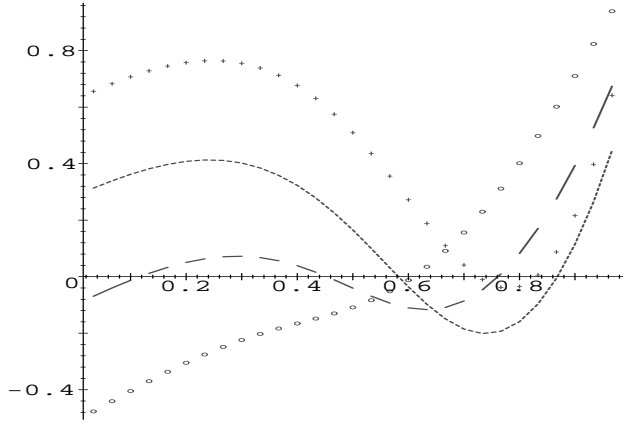}
\caption{We take 3 losing games with $p_0 = p_1
\equiv p$, $q_0 = q_1 \equiv q$, and $r_0 = r_1 \equiv r$, and $p_2 = 0.8(1 - p)^2/((1 -
p)^2 + p^2)$, $q_2 = 0.8( 1- q)^2/((1 - q)^2 + q^2)$, and $r_2 = 0.9( 1-
r)^2/((1 - r)^2 + r^2)$. We set $q = 3/4$, so that $q_2 = 0.08$, and plot
$\ln(c({\bf P}, {\bf Q}, {\bf R})) = d_3 + d_4$ -- compare (\ref{CRIT}), as a
function of $p$ for various values of $r$, circles: $r = 1/2$, $(r_2 = 0.45)$;
 solid: $r = 5/8$, $(r_2 =
0.238235\ldots )$; dashes: $r = 3/4$, $(r_2
= 0.09)$; crosses $r = 7/8$, $(r_2 = 0.018)$. If $\ln(c({\bf P}, {\bf Q}, {\bf
R}))$ is positive (zero) (negative) then the game $({\bf P},{\bf Q}, {\bf R}))$
is winning (fair) (losing).}
\end{figure}

\begin{figure}
\epsfxsize=4in\epsfysize=3in\epsffile{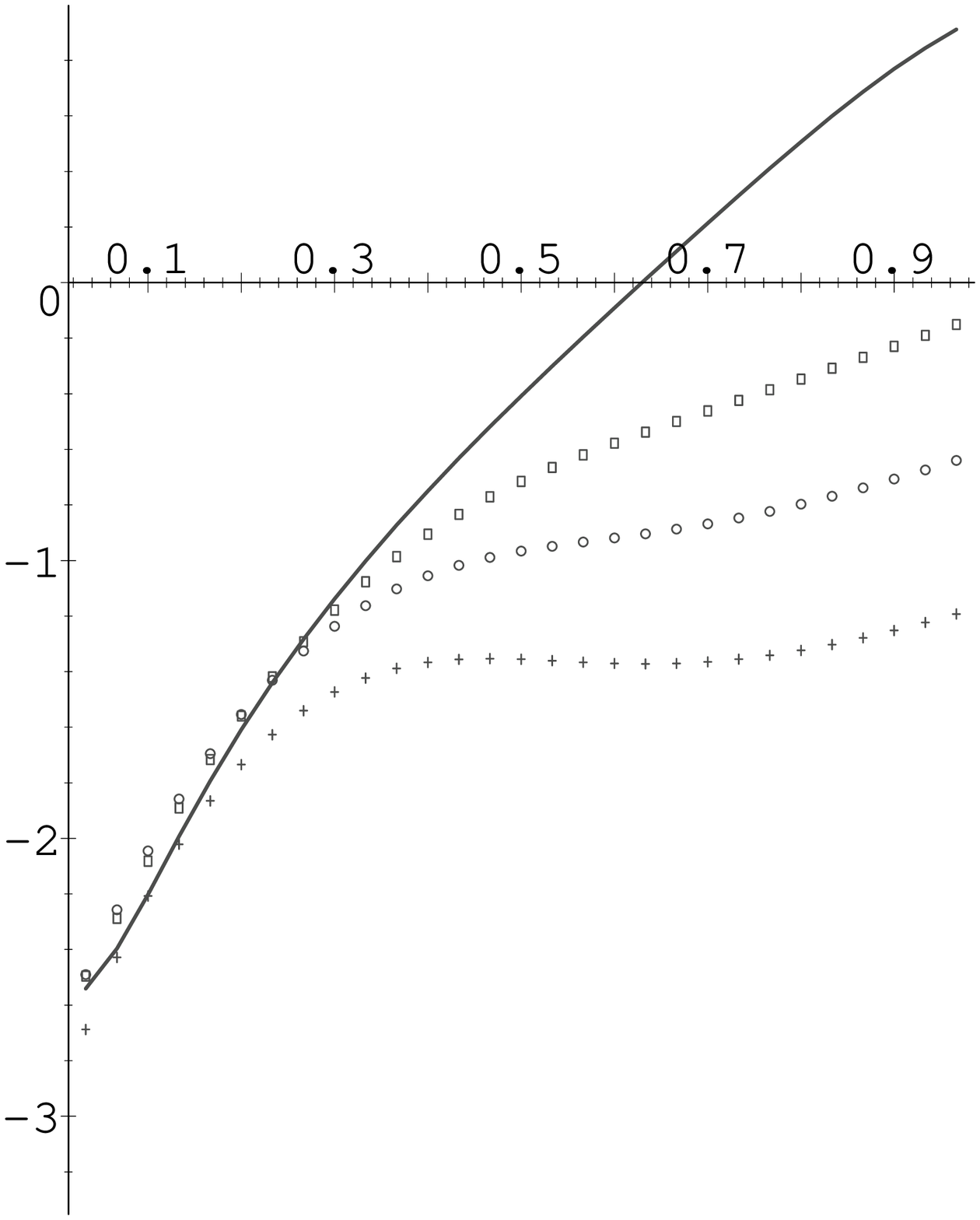}
\caption{We take 3 losing games with $p_0 = p_1
\equiv p$, $q_0 = q_1 \equiv q$, and $r_0 = r_1 \equiv r$, and $p_2 = 0.8(1 - p)^2/((1 -
p)^2 + p^2)$, $q_2 = 0.8( 1- q)^2/((1 - q)^2 + q^2)$, and $r_2 = 0.9( 1-
r)^2/((1 - r)^2 + r^2)$. We set $q = 0.1$ so that $q_2 = 0.8890\ldots$, and
plot $\ln(c({\bf P}, {\bf Q}, {\bf R})) = d_3 + d_4$ -- compare (\ref{CRIT}),
as a function of $p$ for various values of $r$, solid: $r = 3/4$, $(r_2 =
0.9)$; box: $r = 1/2$, $(r_2 = 0.45)$; circles: $r = 3/8$, $(r_2 =
0.66176\ldots)$; crosses $r = 1/4$, $(r_2 = 0.81)$. If $\ln(c({\bf P}, {\bf Q},
{\bf R}))$ is positive (zero) (negative) then the game $({\bf P},{\bf Q}, {\bf
R}))$ is winning (fair) (losing).}
\end{figure}

First we construct a game $({\bf PQR})$ in which the 3--periodic games ${\bf
P}, {\bf Q}$ and ${\bf R}$ are played in a deterministic order. We show that if
all three games are fair then the resulting game may be losing. We also show
that if all three games are winning then the composite game may be losing. If
$X_n$ denotes the winnings of the composite game at time $n$ then its
transition probabilities are given by

\begin{eqnarray}
&&\mbox{if\ } x = 0 \mbox{\ mod\ } 3 \mbox{\ then}\nonumber\\
&&\hspace*{-3ex}P\{X_{n + 3} = x + 3 | X_n = x \} = a_0 \equiv p_0q_1r_2 \nonumber\\
&&\hspace*{-3ex}P\{X_{n + 3} = x - 3 | X_n = x \} = b_0 \equiv (1 - p_0)(1 - q_2)(1 - r_1) \nonumber\\
&&\hspace*{-3ex}P\{X_{n + 3} = x + 1 | X_n = x \} = c_0 \equiv\nonumber\\
&&\hspace*{6ex}p_0( 1- q_2)r_1 + (1 - p_0)q_2(1 - r_0) + p_0(1 - q_1)( 1-
r_0)\nonumber\\
&&\hspace*{-3ex}P\{X_{n + 3} = x - 1 |  X_n = x \} = d_0 \equiv \nonumber\\
 &&\hspace*{6ex}(1 - p_0)(1 - q_2)r_1 + ( 1- p_0)q_2(1 - r_0)  + p_0(1 -
q_1)(1 - r_0)\nonumber
\end{eqnarray}

\begin{eqnarray}
&&\mbox{if\ } x = 1 \mbox{\ mod\ } 3 \mbox{\ then}\nonumber\\
&&\hspace*{-3ex}P\{X_{n + 3} = x + 3 | X_n = x \} = a_1 \equiv p_1q_2r_0\nonumber\\
&&\hspace*{-3ex}P\{X_{n + 3} = x - 3 | X_n = x \} = b_1 \equiv (1 - p_1)(1 - q_0)(1 - r_2)\nonumber\\
&&\hspace*{-3ex}P\{X_{n + 3} = x + 1 | X_n = x \} = c_1
\equiv p_1q_2(1 - r_0) + p_1(1 - q_2)r_1 + (1 - p_1)q_0r_1
\nonumber\\
&&\hspace*{-3ex}P\{X_{n + 3} = x - 1 |  X_n = x \} = d_1 \equiv\nonumber\\
 &&\hspace*{6ex}(1 - p_1)(1 - q_0)r_2 + ( 1- p_1)q_0(1 - r_1) + p_1(1 - q_2)(1
- r_1)\nonumber
\end{eqnarray}

\begin{eqnarray}
&&\mbox{if\ } x = 2 \mbox{\ mod\ } 3 \mbox{\ then}\nonumber\\
&&\hspace*{-3ex}P\{X_{n + 3} = x + 3 | X_n = x \} = p_2q_0r_1\equiv a_2\nonumber\\
&&\hspace*{-3ex}P\{X_{n + 3} = x - 3 | X_n = x \} = (1 - p_2)(1 - q_1)(1 - r_0)\equiv b_2\nonumber\\
&&\hspace*{-3ex}P\{X_{n + 3} = x + 1 | X_n = x \} = p_2q_0(1 - r_1) + p_2(1 - q_0)r_2 + ( 1- p_2)q_1r_2\equiv c_2\nonumber\\
&&\hspace*{-3ex}P\{X_{n + 3} = x - 1 |  X_n = x \} = d_2 \equiv \nonumber\\
 &&\hspace*{6ex} (1 - p_2)(1 - q_1)r_0 + (1 - p_2)q_1(1 - r_2) + p_2(1 - q_0)(1
- r_2).\nonumber
\end{eqnarray}

Hence we define the walk in periodic environment by $Y_n \equiv X_{3n}$, and
determine its transient/recurrent properties by investigating locations of the
eigenvalues of the matrix ${\bf M} = {\bf A}_1{\bf A}_2{\bf A}_0$. Since the
walk $Y_n$ takes on the values $\pm 3, \pm 1$, $R = L = 3$, and each of the
matrices ${\bf A}_i, i = 1, 2, 3$ is a $6 \times 6$ matrix given by

\begin{eqnarray}
{\bf A}_i = \left[\begin{array}{cccccc}
0&-d_i/b_i&1/b_i&-c_i/b_i&0&-a_i/b_i\\
1&0&0&0&0&0\\
0&1&0&0&0&0\\
0&0&1&0&0&0\\
0&0&0&1&0&0\\
0&0&0&0&1&0
\end{array}\right], i = 0, 1, 2.\nonumber
\end{eqnarray}

The examples illustrated in figures 7 and 8
show that the game composed
deterministically from fair (non--symmetric) games can be fair or losing or
winning, depending on parameters of the problem. This form of the paradox is
not observed when composing (deterministically) two or three 2--periodic games.
We observe even more interesting behaviour when composing deterministically unfair
(losing) games. Figures 9 and 10 
show such examples.

By varying just one parameter, it is possible to change the
composed game from losing to winning, change it back to losing, and then back
again to winning -- cf. figure 9. Three state transitions of the composed game have not been
observed before.

\section{Conclusions}

Given the highly non-linear nature of the
recurrence/transience criteria for random walks in periodic environments as a
function of the process parameters, it is not surprising that any scheme that
combines recurrent random walks in periodic environment to produce new random
walks in periodic environments will produce transient processes, except in
trivial cases.

What is surprising to us is that one can find non-trivial cases where random
walks in periodic environments with negative drift can be combined to form
random walks in periodic environments with positive drift, Parrondo's paradox.
We believe that these reversals arise because combining these processes in the
manner described above radically changes the frequency at which transitions
are governed by probability distributions with positive means, and that future
investigations should center on investigating these frequencies.  It is likely
that it will be easier to understand this phenomenon in the case of stochastic
combinations of the processes since these always lead to the study of nearest
neighbour random walks in periodic environments.

When a composite game is constructed from simple 3--periodic games, Parrondo's paradox
is much more interesting and complex than in the case of 2--periodic games.
As shown in figure~9, by varying one parameter
three state transitions of the composed game may be observed. This phenomenon seems to be
similar to changes of the direction of the current in a multiplicative stochastic
ratchet as observed in (K\l osek \& Cox 2000).

As shown in \S3, a game constructed by composing deterministically
2--periodic fair games leads to a pseudo--paradox. In general, a game composed
from an even
number of 2--periodic fair games always exhibits a pseudo--paradox (Key \& K\l osek 2000).
That is, if each simple game is 2--state periodic with an even temporal period, then the
composite game is constructed only from half of the simple games, and the composite
game may be winning, or losing, or fair. Moreover, when
each of the simple games is fair with an odd temporal period, and the composite game is
composed deterministically, then Parrondo's
paradox never occurs. An example in \S3 of a game composed from 2--state periodic and
3--temporal periodic games illustrates this statement.

\begin{acknowledgements}
D.A. would like to acknowledge funding from the Sir Ross and Sir Keith Smith Fund,
GTECH, the Australian Research Council and SA Lotteries. Also thanks are due
to Prof. A.E. Seigman, Stanford University, USA and Prof. J.~Maynard Smith,
University of Sussex, UK for useful communications. The `Boston Interpretation'
originated during a round table discussion with Prof.~H.E.~Stanley's group
at Boston University.
\end{acknowledgements}

\label{lastpage}
\end{document}